\documentclass[11pt,a4paper]{amsart}
\usepackage{amsmath}

\usepackage{amssymb,latexsym}
\usepackage[dvips]{graphicx}

\newtheorem{theorem}{Theorem} [section]

\newtheorem{proposition}{Proposition} [section]

\newtheorem{remark}{Remark}[section]

\title[Extremal metric for $\lambda_1$ on the Klein bottle]
{A unique extremal metric for the least eigenvalue of the Laplacian on the Klein
bottle}

\author[A. El Soufi, H. Giacomini and M. Jazar]{Ahmad El Soufi,
Hector Giacomini and Mustapha Jazar}
\address{Ahmad El Soufi and Hector Giacomini: Universit\'e
Fran\c{c}ois Rabelais de Tours, Laboratoire de Math\'ematiques et
Physique Th\'eorique, UMR-CNRS 6083, Parc de Grandmont, 37200
Tours, France.}
\address{Mustapha Jazar: Lebanese University, Mathematics
Department, P.O. Box 155-012, Beirut, Lebanon.}

\email{ahmad.elsoufi@univ-tours.fr,
hector.giacomini@univ-tours.fr, mjazar@ul.edu.lb}

\keywords{eigenvalue; Laplacian; Klein bottle; extremal metric, Hamiltonian system;
integrable system}

\subjclass[2000]{primary : 58J50; 58E11; 35P15; secondary 37C27}

\begin{document}

\begin{abstract}
We prove the following conjecture recently formulated by Jakobson,
Nadirashvili and Polterovich \cite{JNP}: on the Klein bottle $\mathbb{K}$, the metric of revolution 
$$g_0= {9+ (1+8\cos ^2v)^2\over 1+8\cos ^2v} \left(du^2 + {dv^2\over 1+8\cos ^2v}\right),$$
$0\le u <\frac\pi 2$, $0\le v <\pi$, is the \emph{unique} extremal metric 
of the first eigenvalue of the Laplacian viewed as a functional on the space of all Riemannian metrics of given area.  The proof leads us to study a Hamiltonian dynamical system which turns out to be completely
integrable by quadratures.
\end{abstract}

\maketitle

\section{Introduction and statement of main results}

Among all the possible Riemannian metrics on a compact
differentiable manifold $M$, the most interesting ones are those
which extremize a given Riemannian invariant. In particular, many
recent works have been devoted to the metrics which maximize the
fundamental eigenvalue $\lambda_1(M,g)$ of the Laplace-Beltrami
operator $\Delta_g$ under various constraints (see, for instance,
\cite{ BLY, H, LY, N, P}). Notice that, since $\lambda_1(M,g)$ is
not invariant under scaling ($\lambda_1(M, kg) =
k^{-1}\lambda_1(M,g)$), such constraints are necessary.

In \cite{YY}, Yang and Yau proved that on any compact orientable
surface $M$, the first eigenvalue $\lambda_1(M,g)$ is uniformly
bounded over the set of Riemannian metrics of fixed area. More
precisely, one has, for any Riemannian metric $g$ on $M$,
$$\lambda_1(M,g) A(M,g) \le 8\pi (\hbox {genus} (M) +1),$$
where $A(M,g)$ stands for the Riemannian area of $(M,g)$ (see
\cite{EI1} for an improvement of this upper bound). In the
non-orientable case, the following upper bound follows from Li and
Yau's work \cite {LY}: $\lambda_1(M,g) A(M,g) \le 24\pi (\hbox
{genus} (M) +1)$. On the other hand, if the dimension of $M$ is
greater than 2, then $\lambda_1(M,g)$ is never bounded above over
the set of Riemannian metrics of fixed volume, see \cite{CD}.

Hence, one obtains a relevant topological invariant of surfaces by
setting, for any compact 2-dimensional manifold $M$,
$$ \Lambda (M) = \sup_g \lambda_1(M,g) A(M,g)=\sup_{g\in {\mathcal R}(M)}\lambda_1(M,g),$$
where ${\mathcal R}(M)$ denotes the set of Riemannian metrics of
area 1 on $M$.

On the other hand, in 
spite of the non-differentiability of the functional $g\mapsto\lambda_1(M,g)$ with respect
to metric deformations, a natural notion of extremal (or critical) metric can be
introduced. Indeed, for any smooth
deformation $g_\varepsilon$ of a metric $g$, the
function $\varepsilon \mapsto \lambda_1(M,g_\varepsilon)$ always admits left and right derivatives at $\varepsilon =0$ with
$${d\over d\varepsilon } \lambda_1(M,g_\varepsilon) \Big|_{\varepsilon=0^+}\leq
{d\over d\varepsilon } \lambda_1(M,g_\varepsilon)
\Big|_{\varepsilon=0^-}$$ (see \cite{EI3, EI4} for details). The
metric $g$ is then said to be \emph{extremal} for the functional
$\lambda_1$ under volume preserving deformations if, for any
deformation $g_\varepsilon$ with $g_0=g$ and $vol(M,g_\varepsilon) =
vol(M,g)$, one has
$${d\over d\varepsilon } \lambda_1(M,g_\varepsilon) \Big|_{\varepsilon=0^+}\leq 0\le
{d\over d\varepsilon } \lambda_1(M,g_\varepsilon)
\Big|_{\varepsilon=0^-}.$$
 This last condition can also be formulated as follows:
$$\lambda_1(M,g_\varepsilon)\leq\lambda_1(M,g)+o(\varepsilon)
\mbox{ as } \varepsilon\rightarrow 0.$$

Given a compact surface $M$, the natural questions related to the functional $\lambda_1$ are :

\begin{enumerate}
\item What are the extremal metrics on $M$?
\item Is the supremum $\Lambda (M)$ achieved and, if so, by what
extremal metrics?
\item How does $\Lambda (M)$ depend on (the genus of) $M$?
\end{enumerate}

Concerning the last question, it follows from \cite{CE} that $
\Lambda (M)$ is an increasing function of the genus with a linear
growth rate. Explicit answers to questions (1) and (2) are only
known for the sphere ${\mathbb
S}^2$, the real projective plane ${\mathbb R}P^2$ and the torus
${\mathbb T}^2$. Indeed, the standard metric $g_{{\mathbb
S}^2}$ (resp. $g_{{\mathbb R}P^2} $) is, up to a dilatation, the only extremal metric on ${\mathbb S}^2$ (resp. ${\mathbb R}P^2$) (see \cite{EI2, EI3, MR}) and one has (see \cite{H} and \cite{LY}) 
$$ \Lambda ({\mathbb S}^2) = \lambda_1 ({\mathbb S}^2, g_{{\mathbb S}^2}) A
({\mathbb S}^2, g_{{\mathbb S}^2}) = 8\pi$$ 
and 
$$ \Lambda ({\mathbb R}P^2)
= \lambda_1 ({\mathbb R}P^2, g_{{\mathbb R}P^2}) A({\mathbb R}P^2,
g_{{\mathbb R}P^2} ) = 12\pi.$$
Concerning the torus, the flat metrics $g_{sq}$ and $g_{eq}$ associated respectively with the square lattice ${\mathbb Z}^2$ and the equilateral lattice $ {\mathbb Z}(1,0) \oplus {\mathbb Z}({1\over 2},{\sqrt{3}\over 2}) $ are, up to dilatations, the only extremal metrics on ${\mathbb T}^2$(see \cite{EI3}). Nadirashvili \cite{N} has proved the existence of a regular global maximizer of the functional $g\mapsto\lambda_1({\mathbb T}^2,g)$, which then implies that  
$$ \Lambda ({\mathbb T}^2) = \lambda_1 ({\mathbb T}^2, g_{eq}) A ({\mathbb T}^2, g_{eq}) =
{8\pi^2\over {\sqrt{3}}}.$$ 
However, some steps in Nadirashvili's proof need to be completed as discussed in the recent work of Girouard \cite{G}. The metric $g_{sq}$ corresponds to a saddle point of the functional $\lambda_1 $.
\medskip

{\it What about the Klein bottle $\mathbb{K}$?}

\medskip
 Nadirashvili \cite{N} observed that an extremal metric on $\mathbb{K}$ cannot be a flat metric. Recently, Jakobson, Nadirashvili and
Polterovich \cite{JNP} proved that a metric of revolution 
$$g_0= {9+ (1+8\cos ^2v)^2\over 1+8\cos ^2v} \left(du^2 + {dv^2\over 1+8\cos ^2v}\right),$$
$0\le u <\frac\pi 2$, $0\le v <\pi$, is an extremal metric on $\mathbb{K}$ and 
conjectured that this metric is, up to a dilatation, the unique extremal metric on $\mathbb{K}$.

\medskip

The main purpose of this paper is to prove this conjecture.
Indeed, we will prove the following

\begin{theorem}\label{main}

The Riemannian metric $g_0$ is, up to a dilatation, the unique
extremal metric of the functional $\lambda_1$ under area
preserving deformations of metrics on the Klein bottle $\mathbb{K}$.
\end{theorem}
\begin{remark}\label{rk}
Nadirashvili \cite{N} has given a sketch of proof of the fact that the supremum $\Lambda
(\mathbb{K})$ is necessarily achieved by a regular (real analytic)
Riemannian metric. An immediate consequence of such a result and Theorem \ref{main} would be 
$$\Lambda (\mathbb{K})=
\lambda_1 (\mathbb{K}, g_{0}) A (\mathbb{K}, g_{0})= 12 \pi
E(2\sqrt 2/3)\simeq13.365\,\pi,$$ where $E(2\sqrt 2/3)$ is the complete elliptic integral of the
second kind evaluated at $\frac{2\sqrt2}3$. \end{remark}

It is worth noticing that the metric $g_0$ does not maximize the
systole functional $g\mapsto \hbox{sys}(g)$ (where $\hbox{sys}(g)$
denotes the length of the shortest noncontractible loop) over the
set of metrics of fixed area on the Klein bottle (see \cite {B}),
while on ${\mathbb R}P^2$ and ${\mathbb T}^2$, the functionals
$\lambda_1$ and sys are maximized by the same Riemannian metrics.

The proof of Theorem \ref{main} relies on the characterization of
extremal metrics in terms of minimal
immersions into spheres by the first eigenfunctions. Indeed, 
a metric $g$ is extremal 
for $\lambda_1$ with respect to area preserving deformations
if and only if there exists a family $h_1, \cdots, h_d$
of first eigenfunctions of $\Delta_g$ satisfying $\sum_{i\le d}
dh_i\otimes dh_i =g$ (see \cite{EI3, EI3a}). This last condition actually means that the
map $(h_1, \cdots, h_d):(M,g )\to{\mathbb R}^d$ is an isometric immersion whose image is a minimal immersed submanifold of a sphere. 

As noticed in \cite {JNP}, the surface $(\mathbb{K}, g_{0}) $ is isometrically and minimally immersed in ${\mathbb S}^4$ as  the bipolar
surface of Lawson's minimal torus $\tau_{3,1}$ defined as the
image in ${\mathbb S}^3$ of the map
$$(u, v)\mapsto ( \cos v \exp (3iu),  \sin v \exp (iu)).$$
In fact,  we will prove the following   

\begin{theorem}\label{th1}

The minimal surface $(\mathbb{K}, g_{0}) \hookrightarrow {\mathbb S}^4$ is, up to isometries, the only isometrically and minimally immersed Klein bottle 
into a sphere by its first eigenfunctions.
\end{theorem}

In \cite{EI2}, Ilias and the first author gave a necessary
condition of symmetry for a Riemannian metric to admit isometric
immersions into spheres by the first eigenfunctions. On the Klein
bottle, this condition amounts to the invariance of the metric
under the natural ${\mathbb S}^1$-action on $\mathbb{K}$. Taking
into account this symmetry property and the fact that any metric
$g$ is conformally equivalent to a flat one, for which the
eigenvalues and the eigenfunctions of the Laplacian are explicitly
known, it is of course expected that the existence problem of
minimal isometric immersions into spheres by the first
eigenfunctions reduces to a second order system of ODEs (see
Proposition \ref{prop}). Actually, the substantial part of this
paper is devoted to the study of  the following second order
nonlinear system:
\begin{equation}\label{eq:1}
\left\{\begin{array}{lcl}
\displaystyle{\varphi_1''= (1-2\varphi_1^2 -8\varphi_2^2) \varphi_1}, \\
\\
\displaystyle{\varphi_2'' = (4-2\varphi_1^2 -8\varphi_2^2)
\varphi_2},
\end{array}\right.
\end{equation}
for which we look for periodic solutions satisfying
\begin{equation}\label{eq:3}
\left\{\begin{array}{l}
\varphi_1 \hbox{ is odd and has exactly two zeros in a period,} \\
\varphi_2 \hbox{ is even and positive everywhere;}
\end{array}
\right.
\end{equation}
and the initial conditions
\begin{equation}\label{eq:2}
\left\{\begin{array}{l}
\displaystyle{\varphi_1(0)=\varphi_2'(0)=0\; } \hbox{(from parity conditions (\ref{eq:3})),}\\
\displaystyle{\varphi_2(0)={1\over 2} \varphi_1'(0)=:p \in (0,1 ]
}.
\end{array}
\right.
\end{equation}
 Notice that a similar approach is used in \cite{FP} where the construction of $\mathbb{S}^1$-equivariant
minimal tori in $\mathbb{S}^4$ and $\mathbb{S}^1$-equivariant Willmore tori in $\mathbb{S}^3$ is related to a completely integrable Hamiltonian system.

In \cite{JNP}, Jakobson, Nadirashvili and Polterovich proved that
the initial value $p=\varphi_2(0)=\sqrt{3/8}$ corresponds to a
periodic solution of (\ref{eq:1})-(\ref{eq:2}) satisfying
(\ref{eq:3}). Based on numerical evidence, they conjectured that
this value of $p$ is the only one which corresponds to a periodic
solution satisfying (\ref{eq:3}). As mentioned by them, a
computer-assisted proof of this conjecture is extremely difficult,
due to the lack of stability of the system.

In Section 3, we provide a complete analytic study of System
(\ref{eq:1}). First, we show that this system admits two
independent first integrals (one of them has been already found in
\cite{JNP}). Using a suitable linear change of variables, we show
that the system becomes Hamiltonian and, hence, integrable. The
general theory of integrable Hamiltonian systems tells us that
bounded orbits correspond to periodic or quasi-periodic solutions
(see \cite{Arnold}). However, to distinguish periodic solutions
 from non-periodic ones is not easy in general. Fortunately, our
 first integrals turn out to be
quadratic in the momenta which enables us to apply the classical
Bertrand-Darboux-Whittaker Theorem and, therefore, to completely
decouple the system by means of a parabolic type change of
coordinates $(\varphi_1,\varphi_2)\mapsto(u,v)$.
 We show that, for any
$p\ne\sqrt{3}/2$, the solutions $u$ and $v$ of the decoupled
system are periodic. The couple $(u,v)$ is then periodic if and
only if the periods of $u$ and $v$ are commensurable. We express
the periods of $u$ and $v$ in terms of hyper-elliptic integrals
and study their ratio as a function of $p$. The following fact
(Proposition \ref{pro1}) gives an idea about the complexity of the
situation: there exists a countable dense subset
$\mathcal{P}\subset(0,\sqrt{3}/2)$ such that the solution of
(\ref{eq:1})-(\ref{eq:2}) corresponding to $p\in(0,\sqrt3/2)$ is
periodic if and only if $p\in\mathcal{P}$.

In conclusion, we show that the solution associated with
$p=\sqrt{3/8}$ is the only periodic one to satisfy Condition
(\ref{eq:3}).

\section{Preliminaries: reduction of the problem}

According to \cite{EI3, EI3a}, a necessary and sufficient
condition for a Riemannian metric $g$ on a compact manifold $M$ to be extremal for the functional $\lambda_1$ under
area-preserving metric deformations is that there exists a family
$h_1, \cdots, h_d$ of first eigenfunctions of $\Delta_g$
satisfying
\begin{equation}\label{A4}
\sum_{i\le d} dh_i\otimes dh_i = g,
\end{equation}
which means that the map $h=(h_1, \cdots, h_d)$ is an isometric
immersion from $(M,g)$ to ${\mathbb R}^d$. Since $h_1, \cdots,
h_d$ are eigenfunctions of $\Delta_g$, the image of $h$ is a
minimal immersed submanifold of the Euclidean sphere ${\mathbb
S}^{d-1}\left(\sqrt{2\over \lambda_1 (M,g)}\right)$ of radius
$\sqrt{{2}/{\lambda_1 (M,g)} }$ (Takahashi's theorem \cite{T}). In
particular, we have
\begin{equation}\label{A5}
\sum_{i\le d} h_i^2={2\over \lambda_1 (M,g)}.
\end{equation}

In \cite{EI2}, Ilias and the first author have studied conformal
properties of Riemannian manifolds $(M,g)$ admitting such minimal
isometric immersions into spheres. It follows from their results
that, if $g$ is an extremal metric of $\lambda_1$ under area
preserving deformations, then
\begin{itemize}
\item[(i)] $g$ is, up to a dilatation, the unique extremal metric
in its conformal class, \item[(ii)] $g$ maximizes the restriction
of $\lambda_1$ to the set of metrics conformal to $g$ and having
the same volume, \item[(iii)]the isometry group of $(M,g)$
contains the isometry groups of all the metrics $g'$ conformal to
$g$.
\end{itemize}

For any positive real number $a$, we denote by $\Gamma_{a}$ the
rectangular lattice of ${\mathbb R}^2$ generated by the vectors
$(2\pi, 0)$ and $(0,a)$ and by $ \tilde{g}_a$ the flat Riemannian
metric of the torus $ {\mathbb T}^2_a \simeq {\mathbb R}^2/
\Gamma_{a}$ associated with the rectangular lattice $ \Gamma_{a}$.
The Klein bottle $\mathbb{K}$ is then diffeomorphic to the
quotient of $ {\mathbb T}^2_a $ by the involution $s~: (x,y)
\mapsto (x+\pi, -y)$. We denote by $g_a$ the flat metric induced
on $\mathbb{K}$ by such a diffeomorphism. It is well known that
any Riemannian metric on $\mathbb{K}$ is conformally equivalent to
one of the flat metrics $g_a$.

Let $g=fg_a$ be a Riemannian metric on $\mathbb{K}$. From the
property (iii) above, if $g$ is an extremal metric of $\lambda_1$
under area preserving deformations, then $\hbox {Isom}
(\mathbb{K},g_a)\subset \hbox {Isom} (\mathbb{K},g)$, which
implies that the function $f$ is invariant under the ${\mathbb
S}^1$-action $(x,y)\mapsto (x+t,y)$, $t\in [0,\pi]$, on
$\mathbb{K}$, and then, $f$ (or its lift to ${\mathbb R}^2$) does
not depend on the variable $x$.

\begin{proposition}\label{prop}
Let $a$ be a positive real number and $f$ a positive periodic
function of period $a$. The following assertions are equivalent
\begin{itemize}

\item[(I)] The Riemannian metric $g=f(y)g_a$ on $\mathbb{K}$ is an
extremal metric of the functional $\lambda_1$ under area preserving deformations.
\item[(II)] There exists a homothetic minimal immersion $h=(h_1, \cdots, h_d): (\mathbb{K},g)\to \mathbb{S}^{d-1}$  such that, $\forall i\le d$, $h_i$ is first eigenfunction of $\Delta_g$.

\item[(III)] The function $f$ is proportional to $\varphi_1^2 +4
\varphi_2^2$, where $\varphi_1$ and $\varphi_2$ are two periodic
functions of period $a$ satisfying the following conditions:

\begin{enumerate}
\item[(a)] $(\varphi_1, \varphi_2)$ is a solution of the equations
\begin{equation*}
\left\{\begin{array}{lcl}
\displaystyle{\varphi_1'' = (1-2\varphi_1^2 -8\varphi_2^2) \varphi_1}, \\
\\
\displaystyle{\varphi_2''= (4-2\varphi_1^2
-8\varphi_2^2)\varphi_2};
\end{array}\right.
\end{equation*}
\item[(b)] $\varphi_1$ is odd, $\varphi_2$ is even and
$\varphi'_1(0)=2\varphi_2(0)$; \item[(c)] $\varphi_1$ admits two
zeros in a period and $\varphi_2$ is positive everywhere;
\item[(d)] $\varphi_1^2 +\varphi_2^2 \le 1$ and the equality holds
at exactly two points in a period.

\end{enumerate}
\end{itemize}

\end{proposition}

From the results \cite{EI3, EI3a} mentioned above, it is clear that (I) and (II) are equivalent. Most of the arguments of the proof of ``(II) implies (III)'' can be
found in \cite{N} and \cite {JNP}. For the sake of completeness, we will recall
the main steps. The proof of ``(III) implies (II)'' relies on the fact that
the system (\ref{eq:1}) admits two independent first integrals.

\begin{proof}[Proof of Proposition \ref{prop}] The Laplacian $\Delta_g$ associated with
the Riemannian metric $g=f(y) g_a$ on $\mathbb{K}$ can be
identified with the operator $-{1\over f(y)} \left({\partial_x ^2
}+ {\partial_y ^2 }\right)$ acting on $\Gamma_a$-periodic and
$s$-invariant functions on ${\mathbb R}^2$. Using separation of
variables and Fourier expansions, one can easily show that any
eigenfunction of $\Delta_g$ is a linear combination of functions
of the form $ \varphi_k (y)\cos kx$ and $\varphi_k(y)\sin kx$,
where, $\forall k$, $\varphi_k$ is a periodic function with period
$a$ satisfying $\varphi_k(-y) = (-1)^k \varphi_k(y)$ and
$\varphi_k'' = (k^2-\lambda f)\varphi_k$. Since a first
eigenfunction always admits exactly two nodal domains, the first
eigenspace of $\Delta_g$ is spanned by
$$\left\{\varphi_0(y),\; \varphi_1(y)\cos x,\; \varphi_1(y)\sin x,\;
\varphi_2 (y)\cos 2x,\; \varphi_2(y)\sin 2x \right\},$$ where,
unless they are identically zero, $\varphi_2 $ does not vanish
while $\varphi_0$ and $\varphi_1$ admit exactly two zeros in $[0,
a)$. In particular, the multiplicity of $\lambda_1(\mathbb{K}, g)$
is at most 5.

Let us suppose that $g$ is an extremal metric of $\lambda_1$ under
area preserving deformations and let $h_1, \cdots, h_d$ be a
family of first eigenfunctions satisfying the equations (\ref{A4})
and (\ref{A5}) above. Without loss of generality, we may assume
that $\lambda_1(\mathbb{K}, g) =2$ and that $h_1, \cdots, h_d$ are
linearly independent, which implies that $d\le 5$. Since $h=(h_1,
\cdots, h_d):\mathbb{K}\to {\mathbb S}^{d-1}$ is an immersion, one
has $d\ge 4$. If $d=4$, then using elementary algebraic arguments
like in the proof of Proposition 5 of \cite{MR}, one can see that
there exists an isometry $\rho \in O(4)$ such that $\rho \circ h=(
\varphi_1(y) e^{i x}, \varphi_2(y) e^{2i x})$ with $\varphi_1^2
+\varphi_2^2 =1 $ (eq. (\ref{A5})) and ${\varphi'}_1^2
+{\varphi'}_2^2 =\varphi_1^2 +4\varphi_2^2=f $ (eq. (\ref{A4}))
which is impossible since $\varphi_1^2 +\varphi_2^2 =1 $ implies
that $\varphi_1$ and $\varphi_2$ admit a common critical point.
Therefore, $d=$ multiplicity of $\lambda_1(\mathbb{K}, g) =5$ and
there exists $\rho \in O(5)$ such that $\rho \circ h=(
\varphi_0(y), \varphi_1(y) e^{i x}, \varphi_2 (y) e^{2i x})$, with
$\varphi_0^2+ \varphi_1^2 +\varphi_2^2 =1 $ and
${\varphi'}_0^2+{\varphi'}_1^2 +{\varphi'}_2^2 =\varphi_1^2
+4\varphi_2^2=f $. Since the linear components of $\rho \circ h$
are first eigenfunctions of $(\mathbb{K}, g) $, one should has,
$\forall k=0, 1, 2$, $\varphi_k'' = (k^2-\lambda_1(\mathbb{K}, g)
f)\varphi_k = (k^2-2\varphi_1^2 -8\varphi_2^2) \varphi_k$. Now, it
is immediate to check that one of the couples of functions $
\left(\pm \varphi_1, \pm \varphi_2\right)$ satisfies the
Conditions (a), $\dots$, (d) of the statement. Indeed, the parity
condition $\varphi_k(-y) = (-1)^k \varphi_k(y)$ implies that
$\varphi_1(0)=\varphi'_0(0)=\varphi'_2(0)=0$ and, then,
${\varphi'}^2_1(0)=4\varphi^2_2(0)$. Conditions (c) and (d) follow
from the fact that a first eigenfunction has exactly two nodal
domains in $\mathbb{K}$.

Conversely, let $\varphi_1$ and $\varphi_2$ be two periodic
functions of period $a$ satisfying Conditions (a), $\dots$, (d) of (III)
and consider the Riemannian metric $g=f(y) g_a$ on $\mathbb{K}$, with $f=
\varphi_1^2 +4\varphi_2^2$. We set
$\varphi_0=\sqrt{1- \varphi_1^2 -\varphi_2^2 }$ and define the map
$h:\mathbb{K}\to {\mathbb S}^4$ by $h=( \varphi_0(y), \varphi_1(y)
e^{i x}, \varphi_2 (y) e^{2i x})$. It suffices to check that the
components of $h$ are first eigenfunctions of $\Delta_g$
satisfying (\ref{A4}).

Indeed, in the next section we will see that the second order
differential system satisfied by $\varphi_1$ and $\varphi_2$
(Condition (a)) admits the two following first integrals:
\begin{equation}\label{first0}
\left\{\begin{array}{l}
(\varphi_1^2+4\varphi_2^2)^2-\varphi_1^2-16\varphi_2^2+{\varphi_1'}^2+4{\varphi'_2}^2= C,\\ \\
12\varphi_2^2(\varphi_2^2-1)+3\varphi_1^2\varphi_2^2+
\varphi_2^2{\varphi_1'}^2-2\varphi_1{\varphi_1'}\varphi_2\varphi_2'
+(3+\varphi_1^2){\varphi_2'}^2=C,
\end{array} \right.
\end{equation}
with $C= 4\varphi_2(0)^2(4\varphi_2(0)^2 - 3)$ (note that
Condition (b) implies that $\varphi_1(0)=\varphi'_2(0)=0$).
Differentiating $\varphi_0^2+ \varphi_1^2 +\varphi_2^2 =1 $ and
using the second equation in (\ref{first0}), we get
\begin{eqnarray*}
\varphi_0^2{\varphi_0'}^2&=&\varphi_1^2{\varphi_1'}^2+
\varphi_2^2{\varphi_2'}^2+2\varphi_1\varphi_1'\varphi_2\varphi_2'\\
&=&\varphi_1^2{\varphi_1'}^2+\varphi_2^2{\varphi_2'}^2+
12\varphi_2^2(\varphi_2^2-1)+3\varphi_1^2\varphi_2^2+
\varphi_2^2{\varphi_1'}^2+(3+\varphi_1^2){\varphi_2'}^2-C\\
&=&(\varphi_1^2+\varphi_2^2){\varphi_1'}^2+
(3+\varphi_1^2+\varphi_2^2){\varphi_2'}^2 +
12\varphi_2^2(\varphi_2^2-1)+3\varphi_1^2\varphi_2^2-C\\
&=&(1-\varphi_0^2){\varphi_1'}^2+(4-\varphi_0^2){\varphi_2'}^2
+12\varphi_2^2(\varphi_2^2-1)+3\varphi_1^2\varphi_2^2-C.
\end {eqnarray*}
Therefore
\begin{eqnarray*}
\varphi_0^2\left({\varphi_0'}^2 + {\varphi_1'}^2 +
{\varphi_2'}^2 \right)&=&{\varphi_1'}^2+
4{\varphi_2'}^2+12\varphi_2^2(\varphi_2^2-1)+3\varphi_1^2\varphi_2^2-C\\
&=& \left(1-\varphi_1^2 -\varphi_2^2\right) \left(\varphi_1^2
+4\varphi_2^2\right),
\end {eqnarray*}
where the last equality follows from the first equation of
(\ref{first0}). Hence,
$$|{\partial_y h }|^2= {\varphi_0'}^2 + {\varphi_1'}^2
+{\varphi_2'}^2= \varphi_1^2 +4\varphi_2^2 =
|{\partial_x h}|^2$$
and, since $\partial_x h$ and $\partial_y h$
are orthogonal, the map $h$ is isometric, which means that
Equation (4) is satisfied.

From Condition (a) one has $\varphi_1'' = (1- 2f)\varphi_1$ and
$\varphi_2'' = (4-2f)\varphi_2$, which implies that the functions
$h_1=\varphi_1(y) \cos x$, $h_2=\varphi_1(y) \sin x$,
$h_3=\varphi_2(y) \cos 2x$ and $h_4=\varphi_2(y) \sin 2 x$ are
eigenfunctions of $\Delta_g$ associated with the eigenvalue
$\lambda =2$. Moreover, differentiating twice the identity
$\varphi_0^2+ \varphi_1^2 +\varphi_2^2 =1 $ and using Condition
(a) and the identity ${\varphi_0'}^2 + {\varphi_1'}^2
+{\varphi_2'}^2= \varphi_1^2 +4\varphi_2^2=f$, one obtains after
an elementary computation, $ \varphi_0''=-2f \varphi_0$. Hence,
all the components of $h$ are eigenfunctions of $\Delta_g$
associated with the eigenvalue $\lambda =2$. It remains to prove
that $2$ is the first positive eigenvalue of $\Delta_g$ or,
equivalently, for each $k=0, 1, 2$, the function $\varphi_k$
corresponds to the lowest positive eigenvalue of the
Sturm-Liouville problem $\varphi''=(k^2-\lambda f) \varphi$
subject to the parity condition $\varphi(-y) = (-1)^k \varphi(y)$.
As explained in the proof of Proposition 3.4.1 of \cite{JNP}, this
follows from conditions (c) and (d) giving the number of zeros of
$\varphi_k$, and the special properties of the zero sets of
solutions of Sturm-Liouville equations (oscillation theorems of
Haupt and Sturm).
\end{proof}
\begin{remark}\label{rk1}
Once the initial conditions 
$\varphi'_1(0)=2\varphi_2(0)=p$ and $\varphi_1(0)=\varphi'_2(0)=0$ (since $\varphi_1$ is odd and $\varphi_2$ is even) are fixed, the solution of the system given in assertion (III) of Proposition \ref{prop} is clearly unique. Hence, as we have seen in the proof of this proposition, if a Klein bottle $(\mathbb{K},g=f(y)g_a)$ admits an isometric full minimal immersion $h:(\mathbb{K},g)\to \mathbb{S}^{d-1}$ by the first eigenfunctions, then $d=$the multiplicity of $\lambda_1(\mathbb{K},g)= 5$ and there exists $\rho \in O(5)$ such that
$$\rho \circ h=\left(
\sqrt{1-\varphi_1^2(y)-\varphi_2^2(y)}, \varphi_1(y) e^{i x}, \varphi_2 (y) e^{2i x}\right),$$
where $(\varphi_1,\varphi_2)$ is a unique solution of (III) (with $\varphi'_1(0)=2\varphi_2(0)=\sqrt{f(0)}$ and $\varphi_1(0)=\varphi'_2(0)=0$).
Recall that an immersion $h$ into  $\mathbb{S}^{d-1}$ is said to be \textit{full} if its image does not lie in any hyperplane of  $\mathbb{R}^{d}$ (i.e. its components $h_1,\ldots, h_d$ are linearly independent).

\end{remark}


\section{Study of the dynamical system: proof of results}
According to Proposition \ref{prop}, one needs to deal with the
following system of second order differential equations (Condition
(a) of Prop. \ref{prop})
\begin{equation}\label{s1}\left\{\begin{array}{l}
\varphi_1''=(1-2\varphi_1^2-8\varphi_2^2)\varphi_1,\\
\varphi_2''=(4-2\varphi_1^2-8\varphi_2^2)\varphi_2,
\end{array}\right.\end{equation}
subject to the initial conditions (Condition (b) of Prop.
\ref{prop})
\begin{equation}\label{datas1}\left\{\begin{array}{l}
\varphi_1(0)=0,\; \;\; \varphi_2(0)=p,\\
\varphi_1'(0)=2p,\,  \varphi_2'(0)=0,
\end{array}\right.
\end{equation}
where $p\in (0,1]$ (Condition (d) of Prop. \ref{prop}).

Notice that the system (\ref{s1})-(\ref{datas1}) is invariant
under the transform $$(\varphi_1 (y), \varphi_2(y))\mapsto
(-\varphi_1 (-y), \varphi_2(-y)).$$ Consequently, the solution
$(\varphi_1, \varphi_2)$ of (\ref{s1})-(\ref{datas1})  is such
that $\varphi_1$ is odd and $\varphi_2$ is even.

We are looking for periodic solutions satisfying the following
condition (Condition (c) of Prop. \ref{prop}):
\begin{equation}\label{zeros}
\left\{\begin{array}{l}

\varphi_1 \mbox{ has exactly two zeros in a period,}\\
\varphi_2 \mbox{ is positive everywhere.}
\end{array}
\right.
\end{equation}

Our aim is to prove the following

\begin{theorem}\label{th2}
There exists only one periodic solution of
(\ref{s1})-(\ref{datas1}) satisfying Condition (\ref{zeros}). It
corresponds to the initial value $\varphi_2(0)=p=\sqrt{3/8}$.
\end{theorem}

In fact, this theorem follows from the qualitative behavior of
solutions, in terms of $p$, given in the following

\begin{proposition}\label{pro1} Let $(\varphi_1, \varphi_2)$ be
the solution of (\ref{s1})-(\ref{datas1}).

\begin{enumerate}

\item For all $p\in(0,1]$, $p\ne\sqrt{3}/2$, $(\varphi_1,
\varphi_2)$ is periodic or quasi-periodic.

\item For $p=\frac{\sqrt3}2$, $(\varphi_1, \varphi_2)$ tends to
the origin as $y\to\infty$ (hence, it is neither periodic nor
quasi-periodic).

\item For all $p\in(\sqrt3/2,1]$, $\varphi_2$ vanishes at least
once in each period (of $\varphi_2$). Hence, Condition
(\ref{zeros}) is not satisfied.

\item There exists a countable dense subset
$\mathcal{P}\subset(0,\sqrt{3}/2)$, with $\sqrt{3/8}\in
\mathcal{P}$, such that the solution $(\varphi_1, \varphi_2)$
corresponding to $p\in(0,\sqrt3/2)$ is periodic if and only if
$p\in\mathcal{P}$.

\item For $p=\sqrt{3/8}$, $(\varphi_1, \varphi_2)$ satisfies
(\ref{zeros}) and, for any $p\in \mathcal{P}$, $p\neq\sqrt{3/8}$,
$\varphi_1$ admits at least 6 zeros in a period.

\end{enumerate}
\end{proposition}

Notice that the assertions (2) and (3) of this proposition were
also proved in \cite{JNP} by other methods.

 The first fundamental step in the study of the system above is the existence
of the following two independent first integrals.

\subsection{First integrals}
The functions
\begin{equation}\label{H1H2}
\left\{\begin{array}{l}
H_1(\varphi_1,\varphi_2,\varphi_1',\varphi_2'):=
(\varphi_1^2+4\varphi_2^2)^2-\varphi_1^2-16\varphi_2^2+
(\varphi_1')^2+4(\varphi_2')^2,\\ \\
H_2(\varphi_1,\varphi_2,\varphi_1',\varphi_2'):=
12\varphi_2^2(\varphi_2^2-1)+3\varphi_1^2\varphi_2^2+
\varphi_2^2(\varphi_1')^2\\ \\ \hskip
5cm-2\varphi_1\varphi_1'\varphi_2\varphi_2'
+(3+\varphi_1^2)(\varphi_2')^2,
\end{array} \right.
\end{equation}
are two independent first integrals of (\ref{s1}), i.e. they
satisfy the equation
$$\varphi_1'\frac{\partial H_i}{\partial \varphi_1}+
\varphi_2'\frac{\partial H_i}{\partial \varphi_2}+
\varphi_1''\frac{\partial H_i}{\partial
\varphi_1'}+\varphi_2''\frac{\partial H_i}{\partial
\varphi_2'}\equiv0.$$ The first one, $H_1$, has been obtained by
Jakobson et al. \cite{JNP}. The orbit of a solution of (\ref{s1})
is then contained in an algebraic variety defined by
\begin{equation}
\left\{\begin{array}{l}
H_1(\varphi_1,\varphi_2,\varphi_1',\varphi_2')=K_1,\\ \\
H_2(\varphi_1,\varphi_2,\varphi_1',\varphi_2')=K_2,
\end{array} \right.
\end{equation}
where $K_1$ and $K_2$ are two constants. Taking into account the
initial conditions (\ref{datas1}), one has
$K_1=K_2=-4p^2(3-4p^2)$. In other words, the solution of
(\ref{s1})-(\ref{datas1}) is also solution of
\begin{equation}\label{first}
\left\{\begin{array}{l}
(\varphi_1^2+4\varphi_2^2)^2-\varphi_1^2-16\varphi_2^2+
(\varphi_1')^2+4(\varphi_2')^2+4p^2(3-4p^2)=0,\\ \\
12\varphi_2^2(\varphi_2^2-1)+3\varphi_1^2\varphi_2^2+
\varphi_2^2(\varphi_1')^2-2\varphi_1\varphi_1'\varphi_2\varphi_2'\\
\\ \hskip 4cm+(3+\varphi_1^2)(\varphi_2')^2+4p^2(3-4p^2)=0,
\end{array} \right.
\end{equation}
with the initial conditions
\begin{equation}\label{datafirst}
\left\{\begin{array}{l}
\varphi_1(0)=0,\\ \\
\varphi_2(0)=p.
\end{array} \right.
\end{equation}
Notice that the parameter $p$ appears in both the equations
(\ref{first}) and the initial conditions (\ref{datafirst}). The
system (\ref{first}) gives rise to a ``multi-valued" 2-dimensional
dynamical system in the following way.

\subsection{2-dimensional dynamical systems}\label{2ddynamic}
From (\ref{first}) one can extract explicit expressions of
$\varphi_1'$ and $\varphi_2'$ in terms of $\varphi_1$ and
$\varphi_2$. For instance, eliminating $\varphi_1'$, one obtains
the following fourth degree equation in $\varphi_2'$
\begin{equation}\label{p4}
d_4(\varphi_1,\varphi_2)(\varphi_2')^4-
2d_2(\varphi_1,\varphi_2)(\varphi_2')^2+d_0(\varphi_1,\varphi_2)=0,
\end{equation}
where $d_0$, $d_2$ and $d_4$ are polynomials in $\varphi_1$,
$\varphi_2$ and $p$. The discriminant of (\ref{p4}) is given by
$$\Delta:=-64\varphi_1^2\varphi_2^2w_1w_2w_3,$$
with
$$\begin{array}{l}w_1(\varphi_1,\varphi_2)=\varphi_1^2+\varphi_2^2-1,\\
\\
w_2(\varphi_1,\varphi_2)=p^2\varphi_1^2-(3-4p^2)\varphi_2^2+p^2(3-4p^2),\\
\\
w_3(\varphi_1,\varphi_2)=-(3-4p^2)\varphi_1^2+16p^2\varphi_2^2-4p^2(3-4p^2).
\end{array}$$
It is quite easy to show that, for any $p$, each one of the curves
$(w_i=0)$ contains the orbit of a particular solution of
(\ref{first}). Moreover, the unit circle $(w_1=0)$ represents the
orbit of the solution of (\ref{first}) satisfying the initial
conditions (\ref{datafirst}) with $p=1$. For $p=\sqrt{3/8}$, we
have $w_3\equiv -4w_2$ and the curve $(w_2=0)$ contains the orbit
of the solution of (\ref{first})-(\ref{datafirst}).

These particular algebraic orbits suggest us searching solutions
$(\varphi_1,\varphi_2)$ defined by algebraic relations of the form
$w_4(\varphi_1,\varphi_2)=F(\varphi_1^2,\varphi_2^2)=0$, where $F$
is a polynomial of degree $\le4$. Apart the three quadrics above,
the only additional solution of this type we found is
$$w_4=(\varphi_1^2+4\varphi_2^2)^2-12\varphi_2^2=0.$$
Like $(w_1=0)$, the curve $(w_4=0)$ is independent of $p$ and
represents the orbit of a particular solution of (\ref{first}) for
arbitrary values of $p$. Since
$$w_4(\varphi_1,\varphi_2)=
(\varphi_1^2+4\varphi_2^2-2\sqrt3\,\varphi_2)(\varphi_1^2+4\varphi_2^2+2\sqrt3\,\varphi_2),$$
the set $(w_4=0)$ is the union of two ellipses passing through the
origin, each one being symmetric to the other with respect to the
$\varphi_1$-axis. The upper ellipse
\begin{equation}\label{ellipse}\varphi_1^2+4\varphi_2^2-2\sqrt3\,\varphi_2=0\end{equation}
corresponds to the orbit of the solution of
(\ref{first})-(\ref{datafirst}) associated with
$p=\frac{\sqrt3}2$.

\subsection{Proof of Proposition \ref{pro1}(2): case $\mathbf{ p=\frac{\sqrt3}2}$.}
In this case, the orbit of the solution of
(\ref{first})-(\ref{datafirst}) is given by (\ref{ellipse}). The
only critical point of (\ref{first}) lying on this ellipse is the
origin, which is also a critical point of the system (\ref{s1}).
Therefore, $(\varphi_1(y),\varphi_2(y))$ tends to the origin as
$y$ goes to infinity (see also \cite{JNP}).

\emph{From now on, we will assume that $p\not=\frac{\sqrt3}2$}.

\subsection{A bounded region for the orbit}

The orbit of the solution of (\ref{first})-(\ref{datafirst}) must
lie in the region of the $(\varphi_1,\varphi_2)$-plane where the
discriminant $\Delta$ of  (\ref{p4}) is nonnegative. This region,
$(\Delta\ge0)$, is a bounded domain delimited by the unit circle
$(w_1=0)$ and the quadrics $(w_2=0)$ and $(w_3=0)$. Its shape
depends on the values of $p$.
\begin{itemize}
\item For $p\in(0,\sqrt3/2)$, $(w_2=0)$ and $(w_3=0)$ are
hyperbolas.

\item The case $p=\sqrt{3/8}$ is a special one since then,
$w_3\equiv-4w_2$, and the region $(\Delta\ge0)$ shrinks to the arc
of the hyperbola $(w_2=0)$ lying inside the unit disk.

\item For $p\in(\sqrt3/2,1]$, $w_3$ is positive and $(w_2=0)$ is
an ellipse.

\end{itemize}

From (\ref{p4}) and (\ref{first}) one can express $\varphi_1'$ and
$\varphi_2'$ in terms of $\varphi_1$, $\varphi_2$ and $p$. Thus,
we obtain a multi-valued 2-dimensional dynamical system
parameterized by $p$ with the initial conditions $\varphi_1(0)=0$
and $\varphi_2(0)=p$. However, the dynamics of such a multi-valued
system is very complex to study. Fortunately, as we will see in
the next subsections, the system (\ref{s1}) can be transformed, by
means of a suitable change of variables, into a Hamiltonian
system, completely integrable by quadratures.

\subsection{Hamiltonian dynamical system }

Let us introduce the new variables $q_1$ and $q_2$ defined by
\begin{equation}\label{q12}
q_1:=\frac1{\sqrt2}\,\varphi_1\hskip 1cm,\hskip1cm
q_2:=\sqrt2\,\varphi_2.\end{equation} The system (\ref{s1})
becomes
\begin{equation}\label{s1q}\left\{\begin{array}{l}
q_1''=[1-4(q_1^2+q_2^2)]q_1=-\frac{\partial V}{\partial q_1},\\ \\
q_2''=4[1-q_1^2-q_2^2]q_2=-\frac{\partial V}{\partial q_2},
\end{array}\right.\end{equation}
with
$$V(q_1,q_2):=(q_1^2+q_2^2)^2-\frac12q_1^2-2q_2^2.$$
Therefore, one has a Hamiltonian system with two degrees of
freedom. The Hamiltonian $H$ is given by
$$H(q_1,q_2,q_1',q_2'):=\frac12[(q_1')^2+(q_2')^2]+V(q_1,q_2).$$
This Hamiltonian is a first integral of (\ref{s1q}) (notice that
$H=\frac14 H_1$). A second independent first integral of
(\ref{s1q}) can be obtained from $H_2$. Consequently, the
Hamiltonian system (\ref{s1q}) is integrable and all its bounded
orbits in phase space $(q_1,q_2,q_1', q_2')$ are contained in a
2-dimensional topological torus (see \cite{Arnold}), which means
that the corresponding solutions are periodic or quasi-periodic,
provided that there is no critical point in the closure of the
orbit. However, it is in general difficult to decide whether such
a solution is periodic or not. The corresponding topological torus
obtained from (\ref{first}) is given by:
\begin{equation}\label{firstq12}\left\{\begin{array}{l}
\frac12[(q_1')^2+(q_2')^2]+(q_1^2+q_2^2)^2-\frac12q_1^2-2q_2^2+p^2(3-4p^2)=0,\\
\\
3q_2^2(q_2^2-2)+3q_1^2q_2^2+(q_1')^2q_2^2-2q_1q_1'q_2q_2'\\ \\
\hskip 3cm +\frac12(3+2q_1^2)(q_2')^2+4p^2(3-4p^2)=0.
\end{array}\right.\end{equation}

It is important to notice that the second first integral is also
quadratic in the $q_1'$ and $q_2'$ variables. Indeed, this enables
us to apply the Bertrand-Darboux-Whittaker theorem: \textit{Given
a Hamiltonian system defined by
$$H=\frac12[(q_1)'^2+(q_2')^2]+V(q_1,q_2),$$
the system admits an additional independent first integral,
quadratic in $q_1'$ and $q_2'$, if and only if the system is
separable in cartesian, polar, parabolic, or elliptic-hyperbolic
coordinates.} (see \cite{AP,W} for details).

In our case, an adequate change of variables is a parabolic one,
given by
\begin{equation}\label{parab}
\left\{\begin{array}{l} q_1^2=-\frac23uv,\\ \\
q_2^2=\frac16(3+2u)(3+2v).
\end{array}\right.\end{equation}
Indeed, from (\ref{firstq12}), one obtains after an elementary
computation
\begin{equation}\label{suv}
\left\{\begin{array}{l} (u')^2 =\frac{P(u)}{(u-v)^2},\\ \\
(v')^2 =\frac{P(v)}{(u-v)^2},\end{array} \right.
\end{equation}
where $P(s):=s(1-2s)(3+2s)(2p^2+s)(3-4p^2+2s)$. Observe that
(\ref{suv}) is not completely decoupled yet; this can be done by
means of a suitable change of the independent variable (see
Subsection 3.6). Each one of the quadrics $(w_1=0)$, $(w_2=0)$ and
$(w_3=0)$ is transformed into two parallel lines. Indeed, we have
$w_1(u,v)=-\frac 14(1-2u)(1-2v)$, $w_2(u,v)=-(\frac
32-2p^2+u)(\frac 32-2p^2+v)$ and $w_3(u,v)=4(2p^2+u)(2p^2+v)$.
Also, we have $\Delta=\frac{16}9P(u)P(v)$. Thus, the region
$(\Delta\ge 0)$ is
transformed into the region $(P(u)P(v)\ge 0)$.\\
Observe that the system (\ref{suv}) is symmetric in $u$ and $v$.
As the change of variables (\ref{parab}) is also symmetric in $u$
and $v$, and since $uv$ must be non positive, one can assume,
without loss of generality, that $u\ge0$ and, hence, $-\frac32\le
v\le0$. Now, the condition $(\Delta\ge0)$ implies that $(u,v)\in
I_1\times I_2$, where
\begin{equation}\label{i1}
I_1:=[\alpha_0,1/2]:=\left\{\begin{array}{ll}
\left[0,\frac 12\right]&\mbox{if }p^2<\frac34,\\
\\
\left[2p^2-\frac 32,\frac12\right]&\mbox{if }\frac34<p^2\le1,\\
\end{array}\right.
\end{equation}
and
\begin{equation}\label{i2}
I_2:=[a_0,a_1]:=\left\{\begin{array}{ll}
\left[2p^2-\frac 32,-2p^2\right]&\mbox{if }p^2\le\frac38,\\
\\
\left[-2p^2,2p^2-\frac 32\right]&\mbox{if }\frac38\le p^2<\frac34,\\
\\
\left[-\frac32,0\right]&\mbox{for }\frac34<p^2\le1.\\
\end{array}\right.
\end{equation}
The initial conditions (\ref{datafirst}) become
$$(u(0), v(0))=\left\{\begin{array}{ll}
(0, 2p^2-\frac 32) &\mbox{if }p^2<\frac34,\\
\\
(2p^2-\frac 32, 0) &\mbox{if }\frac34<p^2\le1.\\
\end{array}\right. $$
In all cases, $u(0)$ and $v(0)$ are zeros of $P$. Hence
$$(u'(0), v'(0))=(0,0).$$
The behavior of $(u,v)$ near $y=0$ is then determined by the
acceleration vector
$$(u''(0), v''(0))=\left\{\begin{array}{ll}
(\frac{12p^2}{3-4p^2}, \frac{16p^2(1 - p^2)(3 - 8p^2)}{(3 - 4p^2)}) &\mbox{if }p^2<\frac34,\\
\\
(\frac{16p^2(1
- p^2)(3 - 8p^2)}{(3 - 4p^2)}, \frac{12p^2}{3-4p^2}) &\mbox{if }\frac34<p^2\le1.\\
\end{array}\right. $$

Notice that for $p=\sqrt\frac38$, $v$ is constant, namely
$v(y)=-\frac34$ for all $y$, while for $p=1$, $u$ is constant with
$u(y) = \frac12$ for all $y$.

\subsection{Proof of Proposition \ref{pro1}(1):
Decoupling the system }

In order to completely decouple the previous system we introduce a
change of the independent variable $y\mapsto \tau$ defined by:
$$\frac{d\tau}{dy}=\frac 1{u-v}.$$
Notice that this change of variable is one-to-one since $u-v\ne0$
(indeed, $I_1\cap I_2=\emptyset$). In this new variable, the
system splits into two independent equations:
\begin{equation}\label{su}
(\dot{u})^2 =P(u),
\end{equation}
\begin{equation}\label{sv}
(\dot{v})^2 =P(v),
\end{equation}
where $\dot{u}:=du/d\tau$ and $\dot{v}:=dv/d\tau$. The solution
$\tau\mapsto u(\tau)$ of (\ref{su}) is also a solution of the
second order ODE
\begin{equation}\label{zpp}\ddot{u}=\frac12P'(u),\end{equation} with the
initial conditions $u(0)=\alpha_0$ (see (\ref{i1}) for the
definition of $\alpha_0$) and $\dot{u}(0)=0$, where $P':=dP/du$.
This solution lies on the curve
\begin{equation}\label{ODE}(\dot{u})^2-P(u)=0\end{equation}
in the $(u,\dot{u})$-phase plane of (\ref{zpp}). Since $\alpha_0$
and $\frac12$ are two consecutive zeros of $P$, the equation
(\ref{ODE}) in the region $\alpha_0\le u\le \frac12$ represents a
closed curve. On the other hand, it is easy to check that $P$ and
$P'$ admit no common zero in the interval $[\alpha_0,\frac12]$.
Hence, there exists no critical point for (\ref{zpp}) on the orbit
defined by (\ref{ODE}) inside the region $\alpha_0\le u\le
\frac12$. Consequently, this closed orbit corresponds to a
periodic solution of (\ref{zpp}) and, therefore, $\tau\mapsto
u(\tau)$ oscillates between $\alpha_0$ and $\frac12$.

A similar analysis for $\tau\mapsto v(\tau)$ implies that it is a
periodic solution of
\begin{equation}\label{vpp}\ddot{v}=\frac12P'(v),\end{equation} with the
initial conditions $v(0)=a_0$ (see (\ref{i2}) for the definition
of $a_0$) and $\dot{v}(0)=0$. Consequently, $\tau\mapsto v(\tau)$
oscillates between $a_0$ and $a_1$. This proves Assertion (1) of
Proposition \ref{pro1}(1).

\subsection{Proof of Proposition \ref{pro1}(3): case $\mathbf{p>{\sqrt3\over 2}}$}

We have just seen that $v(\mathbb{R}) = [a_0,a_1]$, with
$a_0=-\frac32$ for $p\in(\frac{\sqrt3}2,1]$ (see (\ref{i2})). This
implies that $q_2$, and then $\varphi_2$, vanishes at least once
in a period (see (\ref{parab})).

\subsection{About the periods of $u$ and $v$: case $\mathbf{p<{\sqrt3\over 2}}$}

Let us denote $\mathcal{T}_u(p)$ the period of $u$. The function
$\tau\mapsto u(\tau)$ oscillates between $\alpha_0=0$ and
$\frac12$ with velocity $(\dot{u})^2 =P(u)\neq 0$ if $u\in
(0,\frac12)$.  Hence, $u(\tau)$ increases from $0$ to $\frac12$
when $\tau$ goes from 0 to $\mathcal{T}_u(p)\over 2$. It follows
that
\begin{equation}\label{Tu}
\mathcal{T}_u(p)=2\int_0^{\frac12}\frac{ds}{\sqrt{P(s)}}.
\end{equation}
Similarly, for $p\ne \sqrt{3/8}$, the period $\mathcal{T}_v(p)$ of
$\tau\mapsto v(\tau)$ is given by
$$
\mathcal{T}_v(p)=2\int_{a_0}^{a_1}\frac{ds}{\sqrt{P(s)}}.
$$
 Setting, for $p\ne \sqrt{3/8}$, $s=(3-8p^2)r-\frac 32+2p^2$,
one can write
\begin{equation}\label{Tv}
\mathcal{T}_v(p)=2\int_{0}^{\frac 12}\frac{dr}{\sqrt{Q(r)}},
\end{equation}
where
$$Q(r):=2r(1-2r)[2p^2+(3-8p^2)r][2-2p^2-(3-8p^2)r][3-4p^2-2(3-8p^2)r].$$
Hence, the functions $\mathcal{T}_u(p)$ and $\mathcal{T}_v(p)$ are
explicitly given by complete hyper-elliptic integrals. Although
the function $\mathcal{T}_v$ is not defined at $p=\sqrt{3/8}$, its
limit exists. Indeed, setting $t=\frac34+r$ and $\alpha:=\frac
34-2p^2$, we get
$$\mathcal{T}_v(p)=2\int_{-\alpha}^{\alpha}\frac{dt}
{\sqrt{(\frac
32-2t)(\frac52-2t)(\frac32+2t)}\sqrt{\alpha^2-t^2}}.$$ As
$\alpha\to0$, we have
$$
\mathcal{T}_v(p)\sim2\int_{-\alpha}^\alpha\frac{4dt}
{3\sqrt{10}\sqrt{\alpha^2-t^2}}.
$$
Thus
$$\lim_{p\to\sqrt{\frac38}}\mathcal{T}_v(p)=\frac{8\pi}{3\sqrt{10}}.$$
On the other hand, we have
$$\mathcal{T}_u(\sqrt{3/8})=\frac45 \Pi(2/5,1/4),$$
where $\Pi$ is the complete elliptic integral of the third kind
given by
$$\Pi(n,m):=\int_0^{\frac\pi2}\frac{d\theta}{(1-n\sin^2\theta)\sqrt{1-m\sin^2\theta}}.$$
Since for $p=\sqrt{3/8}$, $v$ is constant, the
couple $(u,v)$ is periodic of period $\mathcal{T}_u(\sqrt{3/8})$.\\

\noindent\underline{Behavior of ${\mathcal{T}_v-\mathcal{T}_u}$
and $ \mathcal{T}_v/\mathcal{T}_u$ near ${p=0}$}: One has
$$\mathcal{T}_v(p)-\mathcal{T}_u(p)=2\int_0^{\frac 12}\left[
\frac1{\sqrt{Q(s)}}- \frac1{\sqrt{P(s)}}\right]\,ds.$$ The
integral of $ \frac1{\sqrt{P(s)}}$ is singular only at $p=0$. A
direct computation gives
$$\int_0^{\frac 12}\frac{ds}{\sqrt{P(s)}}\sim \int_0^{\frac
12}\frac{ds}{3\sqrt{s^2+2p^2s}}\sim-\frac 23\ln(p).$$ Similarly,
we get
$$\int_0^{\frac 12}\frac{ds}{\sqrt{Q(s)}}\sim -\ln p.$$
In other words $\mathcal{T}_v(p)-\mathcal{T}_u(p)\to+\infty$ as
$p\to0$ while the ratio $ \mathcal{T}_v(p)/\mathcal{T}_u(p)$ goes
to $\frac 32$ (see the figure below).\\

\noindent\underline{Behavior of ${\mathcal{T}_v-\mathcal{T}_u}$
and $ \mathcal{T}_v/\mathcal{T}_u$ near ${p=\sqrt3/2}$}: One has,
as $p\to\sqrt3/2$, $\mathcal{T}_u(p)\sim-\frac23\ln(\sqrt3/2-p)$
and $\mathcal{T}_v(p)\sim-\ln(\sqrt3/2-p)$. Hence,
$\mathcal{T}_v-\mathcal{T}_u\to+\infty$ as $p\to\sqrt3/2$ and
$ \mathcal{T}_v(p)/\mathcal{T}_u(p)$ goes again to $\frac32$.\\

Thus, $\mathcal{T}_v(p)>\mathcal{T}_u(p)$ near $p=0$ and
$p=\sqrt3/2$. Actually, one has
$\mathcal{T}_v(p)>\mathcal{T}_u(p)$ for all $p\in (0,\sqrt3/2)$ as
shown by the graphic representation of $\mathcal{T}_v$ and
$\mathcal{T}_u$ given below.

\begin{center}
\includegraphics[angle=0,width=8cm]{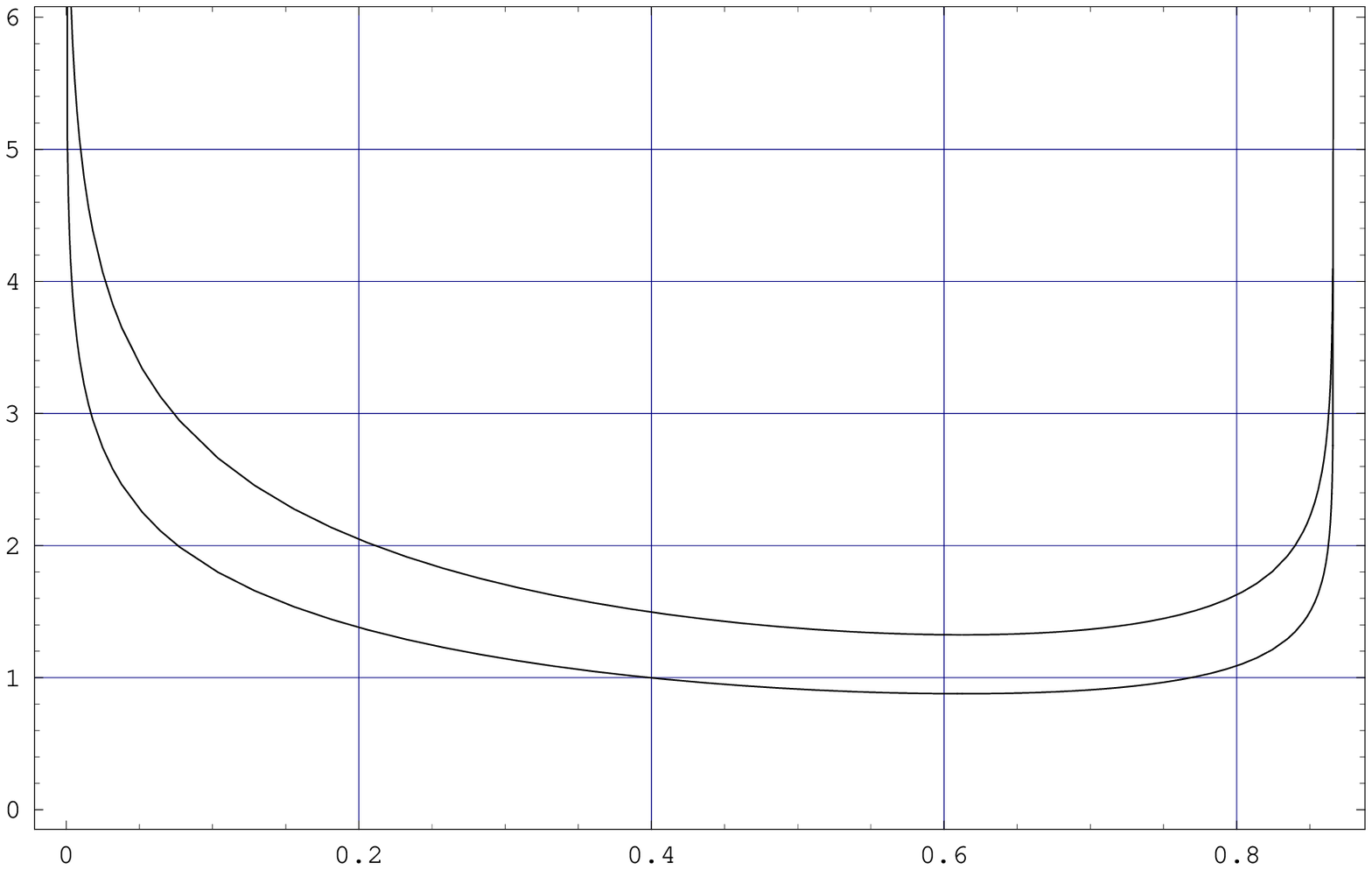}

{\small The functions $p\mapsto \mathcal{T}_v(p)$ (the upper one)
and $p\mapsto\mathcal{T}_u(p)$.}
\end{center}

\subsection{Proof of Proposition \ref{pro1}(4).}
The couple $(u,v)$ is periodic if and only if the ratio
$R(p):=\frac{\mathcal{T}_v(p)}{\mathcal{T}_u(p)}$ is a rational
number. From the previous subsection, $R$ is a nonconstant
continuous function on $(0,\sqrt3/2)$ with $\lim_{p\to
0}R(p)=\lim_{p\to \sqrt3/2}R(p)= 3/2$. The range of $R$ is a
closed interval $[r_1,r_2]\subset[1.480473,1.507784]$ (see the
figure below).

\begin{center}
\includegraphics[angle=0,width=8cm]{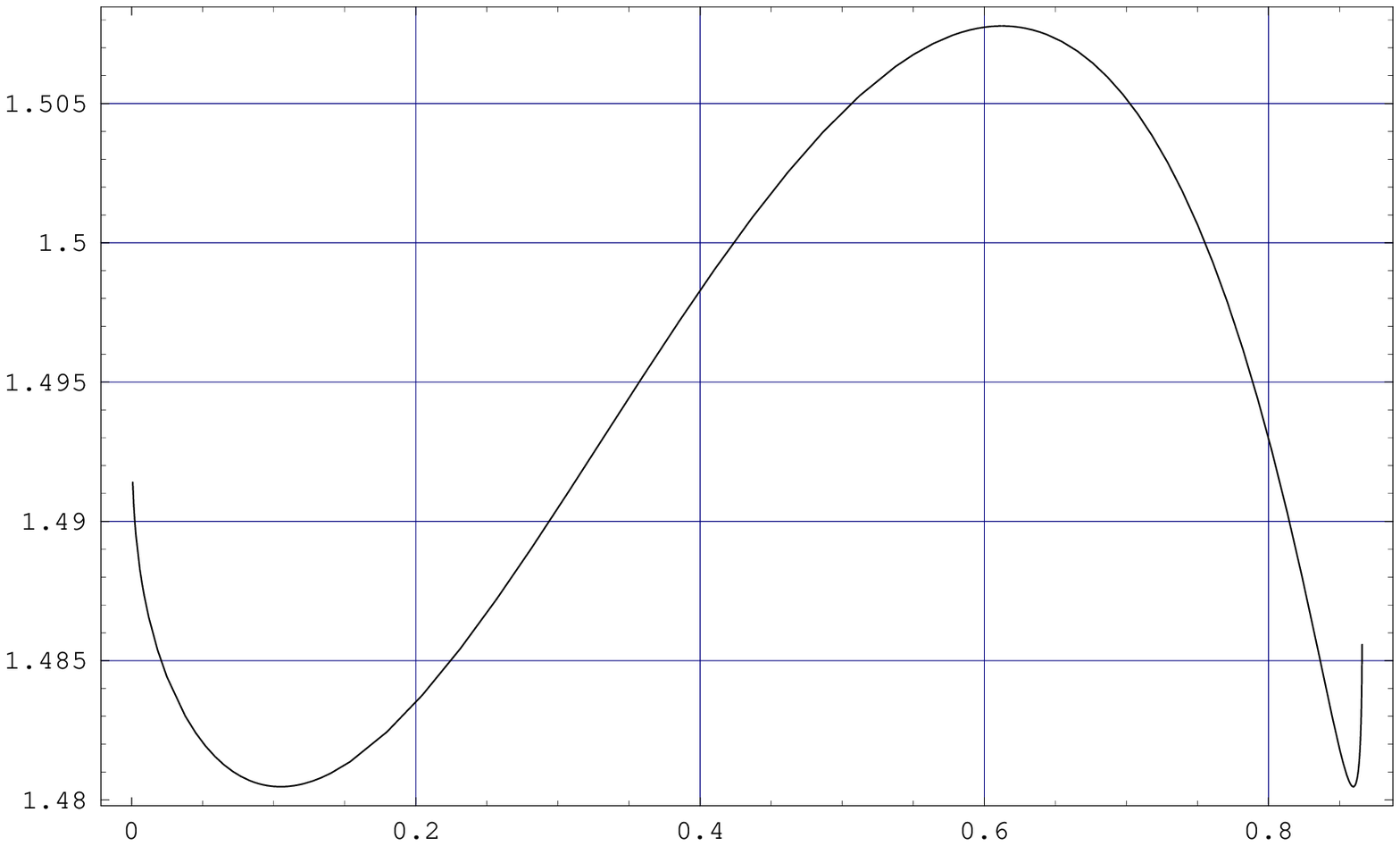}

{\small The ratio $\frac{\mathcal{T}_v}{\mathcal{T}_u}$}
\end{center}

To end the proof of Assertion (4), we only need to define
$\mathcal{P}$ to be the set of $p\in(0,\sqrt3/2)$ such that $R(p)$
is a rational number.

\subsection{Proof of Proposition \ref{pro1}(5).}
In the case $p=\sqrt{3/8}$ we have, $\forall y\in\mathbb{R}$,
$v(y)=-\frac34$ and, then, $\varphi_1^2=u$ and $\varphi_2^2
=\frac18 (3+2u)$. The couple of periodic functions $(\varphi_1,
\varphi_2)=(\sqrt u ,\sqrt{\frac18 (3+2u)})$ on
$[0,\mathcal{T}_u(\sqrt{3/8})]$, such that $\varphi_1$ is odd and
$\varphi_2$ is even, solves the original system and satisfies
Condition (\ref{zeros}).

Let $p\in \mathcal{P}$, $p\neq\sqrt{3/8}$, and let $\frac qm \in
\mathbb{Q}$ be an irreducible fraction, with $q$, $m \in
\mathbb{N}$, such that
$R(p)=\frac{\mathcal{T}_v(p)}{\mathcal{T}_u(p)}=\frac qm$. The
period of the couple $(u,v)$ is given by
$$\mathcal{T}(p)=q\mathcal{T}_u(p)=m\mathcal{T}_v(p).$$
The number of zeros of $u$ in a period, for instance $
[0,\mathcal{T}(p))$, of $(u,v)$ is equal to $q$ times the number
of zeros of $u$ in $[0,\mathcal{T}_u(p))$. As we saw above,
$\forall p\in(0,\sqrt3/2)$, one has $1<R(p)<2$. Hence, $m\ge 2$
and $q>m$, which implies $q\ge3$. Since $u(0)=0$, the number of
zeros of $u$ in a period of $(u,v)$ is at least 3. Since
$\varphi_1$ is odd, the period of $(\varphi_1,\varphi_2)$ is twice
the period of $(u,v)$ (see (\ref{parab}). From  $\varphi_1=
2\sqrt{-\frac13 uv}$ on $[0,\mathcal{T}(p))$, one deduces that
 $\varphi_1$ admits at least 6 zeros in a period of
$(\varphi_1,\varphi_2)$. Notice that the case $p=\sqrt{3/8}$ is
special since, for this value of $p$, $v$ is constant and the
couple $(u,v)$ is periodic whose period is equal to that of $u$.

\subsection{On the shape of solutions}

Although it is not necessary for the proof of our results, one can
obtain as a by product of our study, some properties concerning
the shape of solutions. First, notice that, for
$p\in(0,\frac{\sqrt3}2)$, the 2-dimensional dynamical system
admits four critical points in the region
$(\Delta\ge0)\cap(\varphi_2\ge0)$: $A=(2p/\sqrt3,
\sqrt{1-4p^2/3})$, $B=( \sqrt{1-4p^2/3} , 2p/\sqrt3)$ and their
symmetric with respect to the $\varphi_2$-axis that we denote $A'$
and $B'$. Notice that these critical points are on the boundary of
the region $(\Delta\ge0)$.

\noindent\underline{Non-periodic solutions.} They correspond to
the case where $\mathcal{T}_v(p)/\mathcal{T}_u(p)$ is irrational.
In this case, the orbit $(u,v)$ fill the rectangle $I_1\times
I_2$, and, then, in the $(\varphi_1,\varphi_2)$-plane, the orbit
fill the region $(\Delta\ge0)$.

\noindent\underline{Periodic solutions.}  For $p=\sqrt{3/8}$, the
solution $(\varphi_1,\varphi_2)$ lies on the hyperbola of equation
$\varphi_1^2-4\varphi_2^2+3/2=0,$ oscillating between the points
$A=(-\frac1{\sqrt2},\frac1{\sqrt2})$ and
$A'=(\frac1{\sqrt2},\frac1{\sqrt2})$. For $p\neq\sqrt{3/8}$, let
$q/m=\mathcal{T}_v(p)/\mathcal{T}_u(p)$ be an irreducible
fraction, with $q,m\in\mathbb{N}$, and set
$\mathcal{T}:=q\mathcal{T}_u(p)=m\mathcal{T}_v(p)$. Geometrically,
this means that $u$ makes $q$ round trips in a period while $v$
makes $m$ round trips. We distinguish three cases:
\begin{itemize} \item If $q$ and $m$ are both odd, then
$u(\mathcal{T}/2)=\frac12$ and $v(\mathcal{T}/2)=a_1$. This
corresponds to the point $A$ for $p^2<\frac38$ and to the point
$B$ for $\frac 38<p^2<\frac34$. The orbit is not closed and
$(\varphi_1,\varphi_2)$ oscillates between $A$ and $A'$ or $B$ and
$B'$.

\item If $q$ is odd and $m$ is even, then
$u(\mathcal{T}/2)=\frac12$ and $v(\mathcal{T}/2)=a_0$. This
corresponds to the point $B$ for $p^2<\frac38$ and to the point
$A$ for $\frac 38<p^2<\frac34$. Again, $(\varphi_1,\varphi_2)$
oscillates between $A$ and $A'$ or $B$ and $B'$.

\item If $q$ is even and $m$ is odd, then  $u(\mathcal{T}/2)=0$
and $v(\mathcal{T}/2)=a_1$. This corresponds in the
$(\varphi_1,\varphi_2)$-plane to the point $(0,\frac 34-p^2)$.
This point is the intersection between the quadric $(w_3=0)$ with
the $\varphi_2$-axis. In this case the orbit is closed.
\end{itemize}

\begin{proof}[End of the proof of Theorems \ref{main}, \ref{th1} and \ref{th2}.] 
Theorem \ref{th2} follows directly from Proposition \ref{pro1}.
Let $(\varphi_1,\varphi_2)$ be the periodic solution of
(\ref{s1})-(\ref{datas1}) satisfying (\ref{zeros}) with $\varphi_2(0)=p=\sqrt{3/8}$, and let $a$ be the period of this solution. 
From Proposition \ref{prop}, the Riemannian metric $ \left(\varphi_1^2(y)+\varphi_2^2(y)\right)g_a$ is, up to a dilatation, the only extremal metric of the functional $\lambda_1$ under area preserving deformations.
This proves Theorem \ref{main}. Proposition \ref{prop} also tells us that this metric is the only one to admit a homothetic minimal immersion into a sphere by the first eigenfunctions. Moreover, such an immersion is unique up to isometries and is given by $$\rho \circ h=\left(
\sqrt{1-\varphi_1^2(y)-\varphi_2^2(y)}, \varphi_1(y) e^{i x}, \varphi_2 (y) e^{2i x}\right)$$
for some $\rho \in O(5)$ (see Remark \ref{rk1}).  This proves Theorem \ref{th1}.
\end{proof}

\end{document}